\documentclass[11pt]{article}
\usepackage[T2A]{fontenc}
\usepackage[cp1251]{inputenc}
\usepackage{amsfonts}
\usepackage{eufrak}
\usepackage{amssymb}
\usepackage{amsmath}

\begin{document}

\centerline{\textbf{INDEPENDENT LINEAR STATISTICS ON THE CYLINDERS}}

\bigskip

\centerline{{G.M. FELDMAN and M.V. MYRONYUK}}

\bigskip


\bigskip

 \makebox[20mm]{ }\parbox{125mm}{ \small  Let
either $X=\mathbf{R}\times\mathbf{T}$ or $X=\Sigma_\text{\boldmath $a$}\times\mathbf{T}$,
  where $\mathbf{R}$ is the additive
group of real number,   $\mathbf{T}$ is the cycle group
and   $\Sigma_\text{\boldmath $a$}$ is an $\text{\boldmath $a$}$-adic solenoid. Let
$\alpha_{ij}$, where $i, j=1,2,3,$ be topological automorphisms of
the group $X$. We prove the following analogue of the well-known
Skitovich--Darmois theorem for the group $X$. Let $\xi_j$, where
$j=1, 2, 3$, be independent random variables with values in the
group $X$ and distributions
 $\mu_j$ such that their characteristic functions do not vanish.
If the linear statistics
$L_1=\alpha_{11}\xi_1+\alpha_{12}\xi_2+\alpha_{13}\xi_3$,
$L_2=\alpha_{21}\xi_1+\alpha_{22}\xi_2+\alpha_{23}\xi_3$,
and $L_3=\alpha_{31}\xi_1+\alpha_{32}\xi_2+\alpha_{33}\xi_3$
are independent, then  all $\mu_j$ are  Gaussian distributions.}

 \bigskip


  \bigskip

{\textbf{1. Introduction.}} V.P. Skitovich and G. Darmois proved
the following characterization theorem for a Gaussian distribution
on the real line.

\textbf{Theorem A} (\cite{Skitovich}, \cite{Darmois}, see also
\cite[Ch. 3]{Kag-Lin-Rao}). \textit{Let $\xi_j$, where $j=1,
2,\dots, n,$ and $n\geq 2$, be independent random variables. Let
  $\alpha_j, \beta_j$  be nonzero constants. If the linear statistics
$L_1=\alpha_1\xi_1+\cdots+\alpha_n\xi_n$ and
$L_2=\beta_1\xi_1+\cdots+\beta_n\xi_n$ are independent, then all
random variables $\xi_j$ are Gaussian.}

 A number of works have been devoted to
group analogues of  the Skitovich--Darmois theorem (see e.g.
\cite{Fe7}, \cite{Fe10}, \cite{Fe12}--\cite{Fe3},
\cite{Fe20a}--\cite{FeMy1}, \cite{MiFe-2007}, \cite{NeueScho}, and
also \cite{Fe0}, where one can find additional references). In
particular,
 the following theorem was proved in
 \cite{Fe3}.

\textbf{Theorem B}. {\it Let $X$ be a  second countable locally
compact Abelian group. Assume that  the group  $X$  contains no subgroup
topologically isomorphic to   the circle group. Let $\xi_j$,  where
$j=1,2,\dots,n,$ and $n\geq 2$, be independent random variables
with values in $X$ and   with nonvanishing characteristic
functions. Let $\alpha_j, \ \beta_j$ be topological automorphisms
of
the group  $X$.
 If the linear statistics
$L_1=\alpha_1\xi_1+\cdots+\alpha_n\xi_n$ and
$L_2=\beta_1\xi_1+\cdots+\beta_n\xi_n$ are independent, then all
random variables $\xi_j$ are Gaussian.}

Denote by $\mathbf{T}=\{z\in\mathbf{C}:\ |z|=1\}$ the circle group
(the one-dimensional torus). As have been noted in \cite{Fe-B-1986}, Theorem B fails for the
group $\mathbf{T}$ even for two linear statistics $L_1=\xi_1+\xi_2$
and $L_2=\xi_1-\xi_2$. It follows from this that if a group  $X$
contains a subgroup topologically isomorphic to  $\mathbf{T}$,
then the class of distributions which are characterized by the
independence of linear statistics  $L_1=\alpha_1\xi_1+ \alpha_2\xi_2$ and  $L_2=\beta_1\xi_1+
\beta_2\xi_2$ generally speaking
  contains not only Gaussian
distributions. To put it in another way, if a group  $X$
contains a subgroup topologically isomorphic to  $\mathbf{T}$, the Skitovich--Darmois theorem fails if
we consider two linear statistics of two independent random variables.

The aim of the article is the proof of the following somewhat unexpected result.  We prove that
for some locally compact Abelian groups containing a subgroup topologically isomorphic to
$\mathbf{T}$, the Skitovich--Darmois theorem holds again if
we consider three linear statistics of three independent random variables.

Consider the following general problem.

\textbf{Problem 1.} \textit{Let $X$ be a  second countable locally compact Abelian group,
 $\xi_j$, where
$j=1, 2,\dots, n,$ and $n\geq 2$, be independent random variables
with values in $X$ and     distributions $\mu_j$ such that their
characteristic functions do not vanish. Let $\alpha_{ij}$, where
$i, j =1, 2,\dots,n,$ be topological automorphisms of
  $X$.
 Assume that the linear statistics $L_i=\alpha_{i1}\xi_1+\cdots+\alpha_{in}\xi_n$, where $i=1,
2,\dots, n$, are independent. Describe distributions $\mu_j$.}

It follows from Theorem B that under conditions of Problem 1
if a group $X$  contains no subgroup
topologically isomorphic to   $\mathbf{T}$, then $\mu_j$ are Gaussian.
Therefore we should study only the case when $X$
contains a subgroup topologically isomorphic to $\mathbf{T}$.

As has been proved in \cite{Fe3} the solution of Problem 1 can be
reduced to its solution for a  connected locally compact Abelian
group, one of its characteristic  is its dimension.

If  the dimension of $X$ is equal to
1 and $X$ contains a subgroup topologically isomorphic to $\mathbf{T}$, then
$X\cong \mathbf{T}$. As has bee proved in \cite{Ba-E-Sta}, under the conditions of Problem 1
the
distributions $\mu_j$ are
either Gaussian or convolutions of Gaussian
distributions and signed measures supported in the subgroup
generated by the element of order 2.

 Suppose now that  the dimension of
$X$ is equal to 2 and the group $X$  contains a
subgroup   topologically isomorphic to $\mathbf{T}$. Then the
group
  $X$ is topologically
isomorphic to one of the groups: $\mathbf{R}\times\mathbf{T}$, where $\mathbf{R}$ is the group of real numbers,
$\Sigma_\text{\boldmath $a$}\times\mathbf{T}$, where   $\Sigma_\text{\boldmath $a$}$ is
an  $\text{\boldmath $a$}$-adic solenoid,
and $\mathbf{T}^2$. For
$n=2$ Problem 1 for the groups $\mathbf{R}\times\mathbf{T}$ and
$\Sigma_\text{\boldmath $a$}\times\mathbf{T}$ was solved in
\cite{FeMy1}. As in the case when $X=\mathbf{T}$, the
distributions $\mu_j$ are also either Gaussian or convolutions of Gaussian
distributions and signed measures supported in the subgroup
generated by the element of order 2 in the subgroup $\mathbf{T}$. In the case
when $n=2$ and $X=\mathbf{T}^2$   Problem 1 was not solved
completely, but it is known the complete description of
topological automorphisms $\alpha_j, \beta_j$ of the group  $\mathbf{T}^2$  such that the
independence of the linear statistics $L_1=\alpha_1\xi_1+ \alpha_2\xi_2$ and
$L_2=\beta_1\xi_1+ \beta_2\xi_2$ implies that $\mu_j$ are Gaussian
distributions \cite{MiFe-2007}.

 To the best of our knowledge Problem 1 for $n>2$ has not been
studied yet. Moreover, in all articles devoted to
  the Skitovich--Darmois theorem on groups two linear statistics
of   $n\ge 2$   independent random variables were considered.
But in studying characterization problems
on groups,
in contrast to the classical case, generally speaking, the consideration of
two linear statistics of   $n > 2$   independent random variables
does not make it possible to characterize the Gaussian distributions and their analogues
(see e.g. \cite[\S 14]{Fe0}). As follows from Theorems 1 and 2 proved in the article,
the situation is changed if we consider  $n$ linear statistics of $n$ independent random variables.

In this article we solve Problem 1 completely in the case
when $n=3$ for the group $\mathbf{R}\times\mathbf{T}$, and then, based on this result,
for the group $\Sigma_\text{\boldmath $a$}\times\mathbf{T}$ too.
The answer is somewhat unexpected: the distributions $\mu_j$
must be   Gaussian. To put it in another way, if we consider three linear statistics
of three independent random variables the Skitovich--Darmois theorem holds again
for these groups. We underline that as has  been noted above if we consider two linear statistics
of two independent random variables, then for each of the listed groups the distributions  $\mu_j$
 are  either Gaussian or convolutions of Gaussian
distributions and signed measures, i.e. the Skitovich--Darmois theorem fails.
Moreover, we see that for these groups the class of
distributions which are characterized by the independence of $n$
linear statistics of $n$ independent random variables  depends on $n$.

\bigskip

{\textbf{2. Notation and definitions. The main theorem.}} Let $X$
be a  second countable locally compact Abelian group, $Y$ be its
character group. Let $(x,y)$ be the value of the character $y \in
Y$ at the point $x \in X$. Denote by ${\rm Aut}(X)$ the group of
topological automorphisms of $X$. If $H$ is a closed subgroup of
$Y$, denote by $A(X,H)=\{x \in X: (x,y)=1$ for all $y \in H\}$ its
annihilator. If $X_1$ and $X_2$ are locally compact Abelian
groups, then for any continuous homomorphism $\alpha : X_1 \mapsto
X_2$ define the adjoint homomorphism $\widetilde{\alpha} :Y_2
\mapsto Y_1$ by the formula $(x_1,\widetilde{\alpha}y_2)=(\alpha
x_1, y_2)$ for all $x_1 \in X_1, \ y_2 \in Y_2$. Denote by
$\mathbf{Z}(2)$ the subgroup of $\mathbf{T}$ consisting of
elements $\pm 1$, and denote by $\mathbf{Z}$ the group of
integers. Denote by $I$ the identity automorphism of a group. Put
$Y^{(2)}=\{2y: \ y\in Y\}$. Let $f(y)$ be a function on the group
$Y$, and $h$ be an arbitrary element of $Y$. Denote by $\Delta_h$
the finite difference operator
$$
\Delta_h f(y) = f(y+h) - f(y), \quad y \in Y.
$$
A function $f(y)$ on $Y$ is called a    polynomial  if
$$\Delta_{h}^{n+1}f(y)=0$$ for some $n$ and for all $y, h \in Y$.
The minimal $n$ for which this equality holds is called   the
degree  of the polynomial $f(y)$.
  If $A$ and $B$ are
subsets of $X$, we put $A+B=\{{x\in X:x=a+b, \ a\in A, \ y\in
B}\}$.

Denote by ${M^1}(X)$ the convolution semigroup of probability
distributions on $X$. For ${\mu \in {M^1}(X)}$ we denote by
 $$\widehat
\mu(y) = \int_X (x, y) d\mu(x), \quad y\in Y,
$$
its characteristic function, and define the distribution  $\bar
\mu \in {M}^1(X)$ by the formula $\bar \mu(B) = \mu(-B)$ for
all Borel sets $B$ in $X$. Note that
$\widehat{\bar\mu}(y)=\overline{\widehat\mu(y)}$. Denote by
$\sigma (\mu)$ the support of $\mu\in {M^1}(X)$. If  $H$ is a
closed subgroup of  $Y$ and $\widehat \mu (y)=1$ for $y \in H$,
then $\sigma (\mu) \subset A(X,H)$. We say that a distribution
$\mu\in {M^1}(X)$ is concentrated on a Borel set $A \subset
X$, if $\mu(B)=0$ for any Borel   set $B$ in $X$ such that $B\cap
A=\emptyset$. Denote by $E_x$ the degenerate distribution
concentrated at a point $x \in X$.

A distribution $\gamma \in {M^1}(X)$ is called Gaussian
\cite{PaRaVa}, see also \cite[Ch. IV]{Parthasarathy}, if its
characteristic function is represented in the form
$$
\widehat\gamma(y)=(x,y)\exp\{-\varphi(y)\}, \quad y\in Y,
$$
\noindent where $x\in X$, and $\varphi(y)$ is a continuous
nonnegative function on $Y$ satisfying the equation
\begin{equation}\label{intr2}
    \varphi(u+v)+\varphi(u-v)=2[\varphi(u)+\varphi(v)], \quad u,
    v\in Y.
\end{equation}
Denote by $\Gamma(X)$ the set of Gaussian distributions on $X$. It
should be noted that according to this definition degenerate
distributions are Gaussian.

Let    $X=\mathbf{T}$. Then  $Y\cong\mathbf{Z}.$ In order not to
complicate notation we will assume that $Y=\mathbf{Z}.$
Let $X=\mathbf{R}\times\mathbf{T}$. Denote by $x=(t,z)$, where
$t\in\mathbf{R}$, and $z\in\mathbf{T}$, elements of the group $X$.
If $X=\mathbf{R}\times\mathbf{T}$, then
$Y\cong\mathbf{R}\times\mathbf{Z}$. In order to avoid introducing
new notation we assume that $Y=\mathbf{R}\times\mathbf{Z}$,  and
denote by $y=(s,n)$, where $s\in\mathbf{R}$, and $n\in\mathbf{Z}$,
elements of $Y$. It is easy to verify that every automorphism
$\varepsilon\in{\rm Aut}(\mathbf{R}\times\mathbf{Z})$ is defined
by a matrix $\left(\begin{matrix}a &c\\ 0&p\end{matrix}\right)$, where $a, c\in
\mathbf{R}$,  $a\neq 0$, $p=\pm 1$, and $\varepsilon$ operates on
$Y=\mathbf{R}\times\mathbf{Z}$ as follows
$$
\varepsilon(s,n)=(as+cn,p n),\quad \ s\in\mathbf{R}, \
n\in\mathbf{Z}.
$$
Then the adjoint automorphism
$\delta=\widetilde{\varepsilon}\in{\rm Aut}(X)$ is of the form
$$
\delta(t,z)=(at,e^{ict}z^{p}),\quad \ t\in\mathbf{R}, \
z\in\mathbf{T}.
$$
We identify ${\delta}$ and ${\varepsilon}$ with the
corresponding matrix $\left(\begin{matrix}a &c\\ 0&p\end{matrix}\right)$.

The main result of the article is the following theorem.

\textbf{Theorem 1.} \textit{Let $X=\mathbf{R}\times\mathbf{T}$,
and $\alpha_{ij}\in \mathrm{Aut}(X)$, where $i, j=1,2,3$. Let
$\xi_j$, where $j=1, 2, 3$, be independent random variables with
values in $X$ and    distributions $\mu_j$ such that their
characteristic functions do not vanish. Assume that the linear
statistics $L_1=\alpha_{11}\xi_1+\alpha_{12}\xi_2+\alpha_{13}\xi_3$,
$L_2=\alpha_{21}\xi_1+\alpha_{22}\xi_2+\alpha_{23}\xi_3$, and
$L_3=\alpha_{31}\xi_1+\alpha_{32}\xi_2+\alpha_{33}\xi_3$ are
independent. Then  all $\mu_j$ are either degenerate distributions
or Gaussian distributions  such that their supports are cosets of
a subgroup of $X$ topologically
isomorphic to    $ \mathbf{R}$.}

The proof of Theorem 1 is long enough and it is based on the
series of lemmas. So, we will explain in some words the idea of
the proof. First we prove that under the conditions
of Theorem 1 on an arbitrary locally compact Abelian group  $X$
the convolution
$$\nu=\mu_1*\bar{\mu}_1*\mu_2*\bar{\mu}_2*\mu_3*\bar{\mu}_3$$ is a
Gaussian distribution.  If a group $X$ contains no subgroup
topologically isomorphic to $\mathbf{T}$, then on the group $X$
the Cram\'er theorem about decomposition of a Gaussian
distribution holds \cite {Fe1}, and it implies that the
distributions $\mu_j$ are Gaussian. But on the group
$X=\mathbf{R}\times\mathbf{T}$ a Gaussian distribution may have
non-Gaussian factors. Thus, generally speaking, if   $\nu\in {
M}^1(\mathbf{R}\times\mathbf{T})$, and $\nu$ is a Gaussian
distribution, it does not imply that factors of $\nu$ are also
Gaussian. We prove that under the conditions of Theorem 1 the
Gaussian distribution $\nu$ is not arbitrary, but $\nu$ is
concentrated on a subgroup topologically isomorphic to
$\mathbf{R}$. It follows from this by the classical Cram\'er
theorem about decomposition of a Gaussian distribution that the
distributions $\mu_j$ are also Gaussian.

\bigskip

{\textbf{3. Proof  of lemmas.}} To prove Theorem 1 we need some
lemmas.

\textbf{Lemma 1.} \textit{Let $X$ be a  second countable locally
compact Abelian group, and $\alpha_i, \beta_i \in
\mathrm{Aut}(X)$, where $i = 1, 2$. Let $\xi_j$, where $j=1, 2,
3,$ be independent random variables with values in $X$ and
distributions $\mu_j$. The linear statistics $L_1=\xi_1+\xi_2+\xi_3$,
$L_2=\alpha_1\xi_1+\alpha_2\xi_2+\xi_3$, and
$L_3=\beta_1\xi_1+\beta_2\xi_2+\xi_3$ are independent if and only
if the characteristic functions  $\widehat\mu_j(y)$ satisfy the
equation}
\begin{equation}
\label{u1} \widehat\mu_1(u+\widetilde \alpha_1v+\widetilde
\beta_1w)\widehat\mu_2(u+\widetilde \alpha_2v+ \widetilde
\beta_2w)\widehat\mu_3(u+ v+ w)$$
$$=\widehat\mu_1(u)\widehat\mu_2(u)\widehat\mu_3(u)
\widehat\mu_1(\widetilde \alpha_1v)\widehat\mu_2(\widetilde \alpha_2v)\widehat\mu_3(v)
\widehat\mu_1(\widetilde \beta_1w)\widehat\mu_2(\widetilde
\beta_2w)\widehat\mu_3(w), \quad u, v, w\in Y.
\end{equation}

\textbf{Proof.} Note that if $\xi$ is a random variable with
values in the group $X$  and distribution $\mu$, then
$\widehat\mu(y)=\textbf{ E}[(\xi,y)]$. Obviously, the linear statistics
 $L_1$, $L_2$ and $L_3$ are independent if and only if    the equality
\begin{equation}
\label{l1.1} \textbf{ E}[(\xi_1+\xi_2+\xi_3,
u)(\alpha_1\xi_1+\alpha_2\xi_2+\xi_3,
v)(\beta_1\xi_1+\beta_2\xi_2+\xi_3,w)]$$$$=\textbf{
E}[(\xi_1+\xi_2+\xi_3, u)]\textbf{
E}[(\alpha_1\xi_1+\alpha_2\xi_2+\xi_3, v)]\textbf{
E}[(\beta_1\xi_1+\beta_2\xi_2+\xi_3,w)]
\end{equation}
holds for all $u, v, w \in Y$.  Taking into account that the
random variables  $\xi_j$ are independent, we transform the
left-hand side of equality
 (\ref{l1.1})  as follows
$$ \textbf{ E}[(\xi_1+\xi_2+\xi_3,
u)(\alpha_1\xi_1+\alpha_2\xi_2+\xi_3,
v)(\beta_1\xi_1+\beta_2\xi_2+\xi_3,w)]$$
$$=\textbf{ E}[(\xi_1, u+\widetilde{\alpha}_1 v+ \widetilde{\beta}_1 w)
(\xi_2, u+\widetilde{\alpha}_2 v+ \widetilde{\beta}_2 w) (\xi_3,
u+v+ w)]$$
$$=\textbf{ E}[(\xi_1, u+\widetilde{\alpha}_1 v+ \widetilde{\beta}_1 w)]
\textbf{ E}[(\xi_2, u+\widetilde{\alpha}_2 v+ \widetilde{\beta}_2
w)] \textbf{ E}[(\xi_3, u+v+ w)]$$
$$=\widehat\mu_1(u+\widetilde \alpha_1v+\widetilde
\beta_1w)\widehat\mu_2(u+\widetilde \alpha_2v+ \widetilde
\beta_2w)\widehat\mu_3(u+ v+ w).$$

We transform similarly the right-hand side of equality
(\ref{l1.1})
$$\textbf{ E}[(\xi_1+\xi_2+\xi_3,u)] \textbf{ E}[(\alpha_1\xi_1+\alpha_2\xi_2+\xi_3, v)]
\textbf{E}[(\beta_1\xi_1+\beta_2\xi_2+\xi_3,w)]$$
$$=\textbf{ E}[(\xi_1,u)(\xi_2, u)(\xi_3, u)]
\textbf{ E}[(\xi_1, \widetilde{\alpha}_1
v)(\xi_2,\widetilde{\alpha}_2 v)(\xi_3, v)]
\textbf{E}[(\xi_1,\widetilde{\beta}_1w)(\xi_2,\widetilde{\beta}_2w)(\xi_3,w)]$$
$$=\textbf{ E}[(\xi_1,u)]\textbf{ E}[(\xi_2, u)]\textbf{ E}[(\xi_3,
u)] \textbf{ E}[(\xi_1, \widetilde{\alpha}_1 v)\textbf{
E}[(\xi_2,\widetilde{\alpha}_2 v)$$$$\times\textbf{ E}[(\xi_3, v)]
\textbf{E}[(\xi_1,\widetilde{\beta}_1w)]\textbf{E}[(\xi_2,\widetilde{\beta}_2w)]\textbf{E}[(\xi_3,w)]
$$$$=\widehat\mu_1(u)\widehat\mu_2(u)\widehat\mu_3(u)\widehat\mu_1(\widetilde
\alpha_1v) \widehat\mu_2(\widetilde \alpha_2v)\widehat\mu_3(v)
\widehat\mu_1(\widetilde \beta_1w)\widehat\mu_2(\widetilde
\beta_2w)\widehat\mu_3(w).$$ $\Box$

\textbf{Lemma 2.} \textit{Let $\xi_j$, where $j=1, 2, 3,$ be
independent random variables with values in the group $\mathbf{R}$
and distributions $\mu_j$. Let $a_i, b_i$, where $i=1, 2$, be nonzero
constants. Assume that the linear statistics
$L_1=\xi_1+\xi_2+\xi_3$, $L_2=a_1\xi_1+a_2\xi_2+\xi_3$, and
$L_3=b_1\xi_1+b_2\xi_2+\xi_3$ are independent. Then all $\mu_j$
are either degenerate distributions  or nondegenerate Gaussian
distributions, and the following statements hold:}

 1)
$a_1b_2-a_1a_2b_2-a_1b_1b_2-a_2b_1+a_1a_2b_1+a_2b_1b_2=0;$

  2) \textit{possible combinations of signs of} $a_i, b_i$
\textit{are described in the table}:

\qquad\qquad\qquad\qquad\qquad\qquad\qquad\begin{tabular}{|c|c|c|c|}
  \hline
 $a_1$ & $a_2$ & $b_1$ & $b_2$ \\
  \hline
  + & - & - & +  \\
  + & - & - & -  \\
  - & + & + & -  \\
  - & + & - & -  \\
  - & - & + & -  \\
  - & - & - & +  \\
  \hline
\end{tabular}

  3) $a_1\neq a_2$, $b_1\neq b_2$;

  4) $a_2b_1-a_1b_2\neq 0$;

  5) $\left|%
\begin{array}{cc}
  a_1-1 & a_2-1 \\
  b_1-1 & b_2-1 \\
\end{array}%
\right|\neq 0.$

\textbf{Proof.} By Lemma 1 the characteristic functions
$\widehat\mu_j(y)$ satisfy equation (\ref{u1}), which takes the
form
\begin{equation}
\label{t1.1} \widehat\mu_1(u+a_1v+b_1w)\widehat\mu_2(u+a_2v+
b_2w)\widehat\mu_3(u+ v+ w)$$
$$=\widehat\mu_1(u)\widehat\mu_2(u)\widehat\mu_3(u)\widehat\mu_1(a_1v)
\widehat\mu_2(a_2v)\widehat\mu_3(v)
\widehat\mu_1(b_1w)\widehat\mu_2(b_2w)\widehat\mu_3(w), \quad u,
v, w\in \mathbf{R}.
\end{equation}
By the Skitovich--Darmois theorem the independence of the linear
statistics  $L_1$ and $L_2$ implies that all $\mu_j$ are Gaussian
distributions, i.e. the characteristic functions
 $\widehat\mu_j(s)$
are represented in the form
\begin{equation}\label{f3}
    \widehat\mu_j(s)=\exp\{-\sigma_j s^2 +i \tau_j s\}, \quad s\in \mathbf{R},
\end{equation}
where   $\sigma_j\ge 0$, and $\tau_j\in \mathbf{R}$. Substituting
(\ref{f3}) in (\ref{t1.1}) we get that the following equalities hold:
\begin{equation}\label{t1.2}
    \sigma_1a_1+\sigma_2a_2+\sigma_3=0,
\end{equation}
\begin{equation}\label{t1.3}
    \sigma_1b_1+\sigma_2b_2+\sigma_3=0,
\end{equation}
\begin{equation}\label{t1.4}
    \sigma_1a_1b_1+\sigma_2a_2b_2+\sigma_3=0.
\end{equation}
It follows easily from   (\ref{t1.2})--(\ref{t1.4}) that either
all $\mu_j$ are degenerate distributions or all  $\mu_j$ are
nondegenerate distributions.

 Assume that $\mu_j$
are nondegenerate distributions. Then in
(\ref{t1.2})--(\ref{t1.4}) all $\sigma_j>0$. Prove statements
1)--5).

1) Since in (\ref{t1.2})--(\ref{t1.4}) all $\sigma_j\ne 0$, we have

$$\left|\begin{array}{ccc}
  a_1 & a_2 & 1 \\
  b_1 & b_2 & 1 \\
  a_1b_1 & a_2b_2 & 1 \\
\end{array}\right|=0,$$ i.e. 1) holds.

2) Since all $\sigma_j>0$, statement 2) follows directly from
(\ref{t1.2})--(\ref{t1.4}).

3) Suppose that $a_1=a_2$. Then it follows from 2) that $a_i<0$, and
  (\ref{t1.3}) and (\ref{t1.4}) imply that $(1-a_1)\sigma_3=0$.
Since $a_1<0$, it follows from this that $\sigma_3=0$. Then we get from (\ref{t1.2}) that
 $\sigma_1=\sigma_2=0$. The obtained contradiction proves 3) in the case when $a_1=a_2$.
 In the case when  $b_1=b_2$ we reason similarly.

4) Assume that $a_2b_1-a_1b_2=0$. Then
\begin{equation}\label{f4}
    b_2= a_2b_1/a_1.
\end{equation}
 Taking into account (\ref{f4}), we rewrite (\ref{t1.3}) in the form
$${b_1\over a_1}(\sigma_1a_1+\sigma_2a_2)+\sigma_3=0.
$$
It follows from this and from  (\ref{t1.2}) that
\begin{equation}\label{f5}
    \left(1-{b_1\over a_1}\right)\sigma_3=0.
\end{equation}
 Taking into account 2),  (\ref{f4}) implies that either
 $a_1>0$, \ $a_2<0$, \ $b_1<0$, and $b_2>0$  or
$a_1<0,$ \ $a_2>0$, \ $b_1>0$, and $b_2<0$. Since in the both
cases $a_1b_1<0$, it follows from (\ref{f5}) that $\sigma_3=0$.
Then taking into account that
 $a_ib_i<0$, we get from (\ref{t1.4}) that
$\sigma_1=\sigma_2=0$. The obtained contradiction proves 4).

5) Assume that $\left|%
\begin{array}{cc}
  a_1-1 & a_2-1 \\
  b_1-1 & b_2-1 \\
\end{array}%
\right|= 0$, i.e.
\begin{equation}\label{t1.4.1}
    (a_1-1)(b_2-1)=(a_2-1)(b_1-1).
\end{equation}
Taking into account (\ref{t1.4.1}), we rewrite the left hand-side of equality
 1) in the following way:
$$
a_1b_2-a_1a_2b_2-a_1b_1b_2+a_1a_2b_1b_2-a_2b_1+a_1a_2b_1+a_2b_1b_2-a_1a_2b_1b_2
$$
$$
=a_1b_2(1-a_2-b_1+a_2b_1)-a_2b_1(1-a_1-b_2+a_1b_2)
$$
$$
=a_1b_2(a_2-1)(b_1-1)-a_2b_1(a_1-1)(b_2-1)=(a_2-1)(b_1-1)(a_1b_2-a_2b_1)=0.
$$
Since, according to 4), $a_2b_1-a_1b_2\neq 0$, we have
either $a_2=1$ or $b_1=1$. It follows from 2) that in these cases
$a_1<0$, and $b_2<0$, hence this contradicts to (\ref{t1.4.1}).
$\Box$

\textbf{Lemma 3} (\cite[\S2]{Fe}). \textit{Let $X$ be a
topological Abelian group, and   $G$ be a Borel subgroup of $X$,
$\mu\in
 { M}^1(G)$, and $\mu=\mu_1*\mu_2$. Then the distributions  $\mu_j$
may be replaced by their shifts
 ${\mu}'_j$ such that $\mu=\mu'_1*\mu'_2$,
and  $\mu'_j\in  { M}^1(G)$.}

\textbf{Lemma 4} (\cite{Ba-E-Sta}). \textit{Let $X=\mathbf{T}$, and  $\xi_1$ and
$\xi_2$ be independent random variables with values in the group
$X=\mathbf{T}$ and    distributions $\mu_j$ such that their
characteristic functions do not vanish. Assume that the sum
$\xi_1+\xi_2$ and the difference
 $\xi_1+\xi_2$ are independent. Then
$\mu_j=E_{x_j}*\gamma*\pi_j$, where $x_j \in \mathbf{T}$, $\gamma
\in \Gamma(\mathbf{T})$, $\pi_j$ are signed measures on $\mathbf{Z}(2)$
such that $\pi_1*\pi_2=E_1$. To put it in another way, the characteristic functions
$\widehat\mu_j(n)$ are represented in the form
\begin{equation}\label{f1}
\widehat\mu_1(n)=\exp\{{-\sigma n^2+i\theta_1
n+\kappa(1-(-1)^n)}\}, \quad n\in \mathbf{Z},$$$$
\widehat\mu_2(n)=\exp\{{-\sigma n^2+i\theta_2
n-\kappa(1-(-1)^n)}\}, \quad n\in \mathbf{Z},
\end{equation} where
$\sigma\geq 0$, $\kappa\in\mathbf{R}$, $0\leq \theta_j < 2\pi$,
$j=1, 2$.}

{\bf Lemma 5.} \textit{Let $X=\mathbf{T}$, and  $\xi_j$, where
$j=1, 2, 3,$
 be independent random variables with
values in $X$ and    distributions $\mu_j$ such that their
characteristic functions do not vanish. Let $\alpha_{ij}\in
\mathrm{Aut}(X)$, where $i, j=1, 2, 3$. If the linear statistics
$L_1=\alpha_{11}\xi_1+\alpha_{12}\xi_2+\alpha_{13}\xi_3$,
$L_2=\alpha_{21}\xi_1+\alpha_{22}\xi_2+\alpha_{23}\xi_3$, and
$L_3=\alpha_{31}\xi_1+\alpha_{32}\xi_2+\alpha_{33}\xi_3$ are
independent, then all $\mu_j$ are degenerate distributions.}

\textbf{Proof.} It is obvious that $\mathrm{Aut}(X)=\{I,-I\}$. We
note that, if for some two linear statistics, say for $L_1$ and $L_2$,
we have $L_1=\pm L_2$, then all $\mu_j$ are degenerate
distributions. So, we can assume that $L_i\ne \pm L_j$, where $i,
j=1, 2, 3,$ and $i\ne j$. This easily implies that we can suppose
without loss of generality that
 $L_1=\xi_1+\xi_2+\xi_3$,
$L_2=\xi_1-\xi_2+\xi_3$, and ${L_3=-\xi_1+\xi_2+\xi_3}$.

 Put
$\eta_1=\xi_1+\xi_3$, $\eta_2=\xi_2$. Then $\eta_1$ and $\eta_2$
are independent random variables. The independence of $L_1$ and
$L_2$ implies that the sum $\eta_1+\eta_2$ and the difference
$\eta_1-\eta_2$ are also independent. Applying Lemma 4 to the
random variables $\eta_1$ and $\eta_2$ we obtain that there exists
$\sigma_2\ge 0$ such that
\begin{equation}\label{l4.2}
    |\widehat\mu_1(2n)||\widehat\mu_3(2n)|=|\widehat\mu_2(2n)|=e^{-4\sigma_2 n^2}, \quad
    n\in\mathbf{Z}.
\end{equation}
Reasoning as above we get that the independence of the
linear statistics $L_1$ and $L_3$, and the linear statistics $-L_2$ and
$L_3$ imply that there exist  $\sigma_1\ge 0$ and $\sigma_3\ge 0$
such that
\begin{equation}\label{l4.3}
    |\widehat\mu_2(2n)||\widehat\mu_3(2n)|=|\widehat\mu_1(2n)|=e^{-4\sigma_1 n^2}, \quad
    n\in\mathbf{Z},
\end{equation}
\begin{equation}\label{l4.4}
    |\widehat\mu_1(2n)||\widehat\mu_2(2n)|=|\widehat\mu_3(2n)|=e^{-4\sigma_3 n^2}, \quad
    n\in\mathbf{Z}.
\end{equation}
It follows from (\ref{l4.2})--(\ref{l4.4}) that
$$
\sigma_1+\sigma_3=\sigma_2, \quad \sigma_2+\sigma_3=\sigma_1, \quad \sigma_1+\sigma_2=\sigma_3,
$$
hence
\begin{equation}\label{f12}
    \sigma_1=\sigma_2=\sigma_3=0.
\end{equation}
 Since $\widehat\mu_j(n)$ are the characteristic functions of some distributions,  (\ref{f12})
implies that in (\ref{f1}) $\kappa=0$. Hence, $\mu_j$ are
degenerate distributions. $\Box$

\bigskip

Let $X$ be a locally compact Abelian group, $\alpha_i, \beta_i \in
\mathrm{Aut}(X)$, where $i = 1, 2$. We introduce the notation
$$
L=(\widetilde \alpha_1-I)Y+(\widetilde \beta_1-I)Y, \
M=(\widetilde \alpha_2-I)Y+(\widetilde \beta_2-I)Y, \
N=(\widetilde \alpha_2-\widetilde \alpha_1)Y+(\widetilde
\beta_2-\widetilde \beta_1)Y,
$$
and retain them  in the course of the article.

\textbf{Lemma 6.} \textit{Let $X$ be a  second countable locally
compact Abelian group, and $\alpha_i, \beta_i \in
\mathrm{Aut}(X)$, where $i = 1, 2$. Let $\xi_j$, where $j=1, 2,
3,$ be independent random variables with values in $X$ and
distributions $\mu_j$ such that their characteristic functions do
not vanish.   Assume that the linear statistics
$L_1=\xi_1+\xi_2+\xi_3$, $L_2=\alpha_1\xi_1+\alpha_2\xi_2+\xi_3$,
and $L_3=\beta_1\xi_1+\beta_2\xi_2+\xi_3$ are independent. Put
$\nu_j=\mu_j*\bar\mu_j$. Then $\nu=\nu_1*\nu_2*\nu_3\in
\Gamma(X)$, and the functions $\psi_j(y)=-\log\widehat\nu_j(y)$
satisfy the equations:}
\begin{equation}\label{u13}
\Delta_{h}\Delta_{k}\Delta_{l}\psi_1(y)=0, \quad h, y\in Y, \ k\in
N, \ l\in L,
\end{equation}
\begin{equation}\label{u14}
\Delta_{h}\Delta_{k}\Delta_{l}\psi_2(y)=0, \quad h, y\in Y, \ k\in
N, \ l\in M,
\end{equation}
\begin{equation}\label{u15}
\Delta_{h}\Delta_{k}\Delta_{l}\psi_3(y)=0, \quad h, y\in Y, \ k\in
L, \ l\in M.
\end{equation}

\textbf{Proof.} By Lemma 1 the characteristic functions
 $\widehat\mu_j(y)$  satisfy equation (\ref{u1}). Note that $\widehat\nu_j(y)=|\widehat\mu_j(y)|^2>0$,
$y\in Y$, and
 the characteristic functions
 $\widehat\nu_j(y)$  also satisfy
equation (\ref{u1}). It follows from equation (\ref{u1}) that
\begin{equation}\label{u4}
\psi_1(u+\widetilde \alpha_1v+\widetilde
\beta_1w)+\psi_2(u+\widetilde \alpha_2v+\widetilde
\beta_2w)+\psi_3(u+v+w)$$$$=A(u)+B(v)+C(w), \quad u, v, w\in Y,
\end{equation}
where $A(y)=\psi_1(y)+\psi_2(y)+\psi_3(y)$, \
$B(y)=\psi_1(\widetilde \alpha_1y)+\psi_2(\widetilde
\alpha_2y)+\psi_3(y)$, $C(y)=\psi_1(\widetilde
\beta_1y)+\psi_2(\widetilde \beta_2y)+\psi_3(y)$.

We use the finite difference method to solve equation (\ref{u4}).
Let $k_1$ and $m_1$ be   arbitrary elements of the group $Y$. Set
$h_1=-k_1-m_1$, then $h_1+k_1+m_1=0$. Substitute $u+h_1$ for $u$,
$v+k_1$ for $v$, and  $w+m_1$ for $w$ in equation (\ref{u4}).
Subtracting equation (\ref{u4}) from the resulting equation we
obtain
\begin{equation}\label{u5}
\Delta_{l_{11}}\psi_1(u+\widetilde \alpha_1v+\widetilde
\beta_1w)+\Delta_{l_{12}}\psi_2(u+\widetilde \alpha_2v+\widetilde
\beta_2w)$$$$=\Delta_{h_1}A(u)+\Delta_{k_1}B(v)+\Delta_{m_1}C(w),
\quad u, v, w\in Y,
\end{equation}
where ${l_{11}}=(\widetilde \alpha_1-I)k_1+(\widetilde
\beta_1-I)m_1$, and ${l_{12}}=(\widetilde
\alpha_2-I)k_1+(\widetilde \beta_2-I)m_1$. Let $k_2$ and $m_2$ be
arbitrary elements of the group $Y$. Set
  $h_2=-\widetilde \alpha_2k_2-\widetilde
\beta_2m_2$. Then $h_2+\widetilde \alpha_2k_2+\widetilde
\beta_2m_2=0$. Substitute $u+h_2$ for $u$, $v+k_2$ for $v$, and
$w+m_2$ for $w$ in equation (\ref{u5}). Subtracting equation
(\ref{u5}) from the resulting equation we find
\begin{equation}\label{u6}
\texttt{}\Delta_{l_{21}}\Delta_{l_{11}}\psi_1(u+\widetilde
\alpha_1v+\widetilde \beta_1w)=\Delta_{h_2}\Delta_{h_1}A(u)+
\Delta_{k_2}\Delta_{k_1}B(v)+\Delta_{m_2}\Delta_{m_1}C(w), \quad
u, v, w\in Y,
\end{equation}
where ${l_{21}}=(\widetilde \alpha_2-\widetilde
\alpha_1)k_2+(\widetilde \beta_2-\widetilde \beta_1)m_2$. Let
$k_3$ and $m_3$ be arbitrary elements of the group $Y$. Set
  $h_3=-\widetilde \alpha_1k_3-\widetilde \beta_1m_3$.
  Then $h_3+\widetilde \alpha_1k_3+\widetilde \beta_1m_3=0$.
  Substitute $u+h_3$ for $u$,
$v+k_3$ for $v$, and  $w+m_3$ for $w$ in equation (\ref{u6}).
Subtracting equation (\ref{u6}) from the resulting equation we get
\begin{equation}\label{u7}
\Delta_{h_3}\Delta_{h_2}\Delta_{h_1}A(u)+
\Delta_{k_3}\Delta_{k_2}\Delta_{k_1}B(v)+\Delta_{m_3}\Delta_{m_2}\Delta_{m_1}C(w)=0,
\quad u, v, w\in Y.
\end{equation}
Let $h_4$ be an arbitrary element of the group $Y$.  Substitute
$u+h_4$ for $u$ in equation (\ref{u7}). Subtracting equation
(\ref{u7}) from the resulting equation we obtain
\begin{equation}\label{u8}
\Delta_{h_4}\Delta_{h_3}\Delta_{h_2}\Delta_{h_1}A(u)=0, \quad h_j,
u\in Y.
\end{equation}
Putting in (\ref{u8}) $h_1=h_2=h_3=h_4=h$, we receive
\begin{equation}\label{u9}
\Delta^4_{h}A(u)=0, \quad h, u\in Y.
\end{equation}
Thus, $A(y)$ is a polynomial of degree $\le 3$. Then there exist
symmetric $k$-additive functions $g_k(y_1, y_2, \dots, y_k)$,
where $k=1, 2, 3,$ such that

$$ A(y)= g_3^*(y)+g_2^*(y)+g_1^*(y)+g_0,\quad y\in Y,$$
where $g_k^*(y)=g_k(y,  \dots, y)$, and $g_0=const$ \cite{Djo}.
It is obvious that $g_k^*(-y)=(-1)^k g_k^*(y)$. Since $A(-y)=A(y)$
and $A(0)=0$, we have $A(y)=g_2^*(y)$. It is not difficult to
verify that the function $g_2^*(y)$ satisfies equation
(\ref{intr2}). Hence, $\nu \in \Gamma(X)$.

Since
$A(y)=g_2^*(y)=g_2(y, y)$, it is easily seen that the function
$A(y)$ satisfies the equation
\begin{equation}\label{u11}
\Delta_{h_3}\Delta_{h_2}\Delta_{h_1}A(u)=0, \quad h_j, u\in Y.
\end{equation}
Reasoning as above we get that the functions $B(v)$ and $C(w)$
also satisfy equation (\ref{u11}).

Let $h$ be an arbitrary element of the group $Y$.  Substitute
$u+h$ for $u$ in equation (\ref{u6}). Subtracting equation
(\ref{u6}) from the resulting equation and taking into account
(\ref{u11}) we find
\begin{equation}\label{u12}
\Delta_{h}\Delta_{l_{21}}\Delta_{l_{11}}\psi_1(u+\widetilde
\alpha_1v+\widetilde \beta_1w)=0.
\end{equation}
Putting here $v=w=0$ and considering expressions for $l_{11}$ and
$l_{21}$, we obtain that the function $\psi_1(y)$ satisfies
equation (\ref{u13}). Reasoning as above we make sure that the
functions $\psi_2(y)$ and $\psi_3(y)$ satisfy equations
(\ref{u14})  and (\ref{u15}) respectively. $\Box$

\textbf{Lemma 7.} \textit{Let $X=\mathbf{R}\times\mathbf{T}$, and
$\alpha_i = \left(\begin{matrix}a_i&c_i\\ 0&p_i\end{matrix}\right)$,
 $\beta_i =
\left(\begin{matrix}b_i&d_i\\ 0&q_i\end{matrix}\right),$ where $i=1, 2, $ are
topological automorphisms of the group  $X$.  Let $\xi_j$, where
$j=1, 2, 3,$ be independent random variables with values in $X$
and    distributions $\mu_j$.   Assume that the linear statistics
$L_1=\xi_1+\xi_2+\xi_3$, $L_2=\alpha_1\xi_1+\alpha_2\xi_2+\xi_3$,
and $L_3=\beta_1\xi_1+\beta_2\xi_2+\xi_3$ are independent. Then
either all $\mu_j$ are  degenerate distributions  or all $\mu_j$ are
nondegenerate
distributions, and there are the following possibilities for the
subgroups $L$, $M$, and $N$:}

$$\begin{array}{l}{\bf 1.} \ L=\mathbf{R},  \ M=\mathbf{R},  \ N=\mathbf{R}; \\

{\bf 2.} \ L=\mathbf{R}, \ M=Y^{(2)}, \ N=Y^{(2)};\\

{\bf 3.} \ L=Y^{(2)}, \  M=\mathbf{R}, \ N=Y^{(2)}; \\

{\bf 4.} \ L=Y^{(2)}, \ M=Y^{(2)}, \ N=\mathbf{R};\\

{\bf 5.} \ L=Y^{(2)}, \ M=Y^{(2)}, \ N=Y^{(2)}.
\end{array}
$$

\textbf{Proof.} By Lemma 1 the characteristic functions
$\widehat\mu_j(y)$ satisfy equation (\ref{u1}). Putting in
(\ref{u1}) $u=(s_1, 0), \ v=(s_2, 0), \ w=(s_3, 0)$, we obtain
\begin{equation}\label{u2}
 \widehat\mu_1(s_1+a_1s_2+b_1s_3, 0)\widehat\mu_2(s_1+a_2s_2+b_2s_3, 0)
 \widehat\mu_3(s_1+s_2+s_3, 0)=\widehat\mu_1(s_1, 0)\widehat\mu_2(s_1, 0)
 \widehat\mu_3(s_1, 0)$$$$\times\widehat\mu_1(a_1s_2, 0)\widehat\mu_2(a_2s_2, 0)\widehat\mu_3( s_2, 0)
\widehat\mu_1(b_1s_3, 0)\widehat\mu_2(b_2s_3, 0)\widehat\mu_3(
s_3, 0), \quad s_1, s_2, s_3\in \mathbf{R}.
\end{equation}
 Set
$\nu_j=\mu_j*\bar{\mu}_j$. Then
$\widehat\nu_j(y)=|\widehat\mu_j(y)|^2>0$,  $y\in Y$. Obviously,
the functions $\widehat\nu_j(y)$ also satisfy equation (\ref{u1}).
Hence, they satisfy equation (\ref{u2}). It follows from
(\ref{u2}) and Lemma 2 for the group $\mathbf{R}$ that if at least
one of the functions $\widehat\nu_j(s, 0)$ is the characteristic
function of a nondegenerate distribution, then all
$\widehat\nu_j(s, 0)$ are the characteristic function of
nondegenerate Gaussian distributions. To put it in another way, either
$\widehat\nu_j(s,0)=1$, \ $s\in \mathbf{R}$ or
$\widehat\nu_j(s,0)=e^{-\sigma_j s^2}$,   \ $s\in \mathbf{R},$
where $\sigma_j> 0$, and $j=1, 2, 3$.

Assume that $\widehat\nu_j(s,0)=1$, \ $s\in \mathbf{R},$ where
$j=1, 2, 3$.
 This implies the inclusions  ${\sigma(\nu_j)\subset A(X, \mathbf{R})}=\mathbf{T},$ where $j=1, 2, 3$.
 Hence, by Lemma 3   the distributions  $\mu_j$
may be replaced by their shifts
 ${\mu}'_j$ such that  $\sigma(\mu_j')\subset\mathbf{T},$ where $j=1, 2, 3$.
 Since the subgroup $\mathbf{T}$ is invariant with respect to any
 topological automorphism of the group $X$, by Lemma 5 all  $\mu_j$ are degenerate distributions.

Suppose now that  $\widehat\nu_j(s,0)=e^{-\sigma_j s^2}$,  where
 $s\in \mathbf{R},$ $\sigma_j> 0,$  and $j=1, 2, 3$. Then statement 2) of Lemma 2
implies that the following conditions hold:
$$
\begin{array}{l}(i) \ \ \ \mbox{either } a_1\ne 1  \ \mbox{or } b_1\ne 1;\\
(ii) \ \ \mbox{either } a_2\ne 1  \ \mbox{or }   b_2\ne 1.\\
\end{array}
$$

Note that it follows from  $(i)$ that there are only two
possibilities for $L$: either $L=\mathbf{R}$  or $L=Y^{(2)}$.
Reasoning as above we get from $(ii)$ that there are only two
possibilities for $M$: either $ M=\mathbf{R}$  or $M=Y^{(2)}$.
Obviously, each of these possibilities can be realized.

Let  $L=\mathbf{R}$ и  $ M=\mathbf{R}$. Then $p_i=q_i=1,$ where
$i=1, 2.$ Hence, $p_2-p_1=0$, and $q_2-q_1=0$. Taking into account
statement 3) of Lemma 2, this implies that for $N$ there exists
the only possibility  $N=\mathbf{R}$.

Let  $L=\mathbf{R}$ and  $M=Y^{(2)}$. Then $p_1=q_1=1$ and either
$p_2=-1$  or $q_2=-1$. Hence, either $p_2-p_1=-2$  or
$q_2-q_1=-2$. Taking into account statement 3) of Lemma 2, this
implies that for $N$ there exists the only possibility
$N=Y^{(2)}$.

Let  $L=Y^{(2)}$ and  $ M=\mathbf{R}$. Then either $p_1=-1$  or
$q_1=-1$, and $p_2=q_2=1$. Hence, either $p_2-p_1=2$  or
$q_2-q_1=2$. Taking into account statement 3) of Lemma 2, this
implies that for $N$ there exists also the only possibility
$N=Y^{(2)}$.

Let  $L=Y^{(2)}$ and $ M=Y^{(2)}$. Then either $p_1=-1$  or
$q_1=-1$,  and either $p_2=-1$ or $q_2=-1$. If $p_1=p_2$ and
$q_1=q_2$, then by statement 3) of Lemma 2,
 $N=\mathbf{R}$. If either  $p_1\ne p_2$  or  $q_1\ne q_2$, then by statement 3) of Lemma
 2,
 $N=Y^{(2)}$. $\Box$

\textbf{Lemma 8.} \textit{Let $Y=\mathbf{R}\times\mathbf{Z}$, and
a continuous function $f(y)$, $y=(s, n)\in Y,$ satisfies the
equation
\begin{equation}\label{u16}
 \Delta^2_k\Delta_hf(y)=0, \quad k\in  \mathbf{R}, \ h, y\in Y.
\end{equation}
Assume also that $f(-y)=f(y)$, $y\in Y$. Then
\begin{equation}\label{u17}
 f(s, n)=\sigma s^2+\kappa(n)s+\lambda(n), \quad s\in\mathbf{R}, \
n\in\mathbf{Z},
\end{equation}
where  $\kappa(-n)=-\kappa(n)$ and $\lambda(-n)=\lambda(n)$, $n\in
\mathbf{Z}$.}

\textbf{Proof.} Put $f_n(s)=f(s, n),$ where $s\in\mathbf{R}, \
n\in\mathbf{Z}$. Putting in  (\ref{u16}) $y=(s, n)$, $k=h=(t, 0)$,
we get that the function $f_n(s)$ for any fixed  $n$ satisfies the
equation
\begin{equation}\label{u18a}
 \Delta^3_tf_n(s)=0, \quad t, s\in  \mathbf{R}.
\end{equation}
It follows from  (\ref{u18a}) that the function $f_n(s)$ is of the
form
\begin{equation}\label{u18}
f_n(s)=\sigma(n) s^2+\kappa(n)s+\lambda(n), \quad s\in\mathbf{R}.
\end{equation}
Substitute representation (\ref{u18}) for the function $f_n(s)$
into (\ref{u16}), supposing that $h=(0, m)$, $k=(t, 0)$. We get
$$\Delta^2_k\Delta_hf_n(s)=\Delta_h\Delta^2_tf_n(s)=
\Delta_h\sigma(n)\Delta^2_ts^2=2t^2\Delta_m\sigma(n)=0, \quad
t\in\mathbf{R}, \ m\in\mathbf{Z}.$$ This implies that
$\Delta_m\sigma(n)=0$, i.e. $\sigma(n)=\sigma=const$.

Since $f(-y)=f(y)$, \ $y\in Y$, it is obvious that
 $\kappa(-n)=-\kappa(n)$ and $\lambda(-n)=\lambda(n)$,\  $n\in \mathbf{Z}$. It should be noted
 that any function $f(s, n)$ of the form (\ref{u17}) satisfies equation (\ref{u16}). $\Box$

\bigskip

{\textbf{4. Proof of Theorem 1.}} Let $X$ be an arbitrary second
countable locally compact Abelian group. First we note that if
$\mu$ is the distribution of a random variable $\xi$ with values
in the group   $X$ and $\alpha \in {\rm Aut}(X)$, then the
characteristic function of a random variable $\alpha \xi$ is
$\widehat \mu(\widetilde \alpha y)$. This implies that $\mu\in
\Gamma(X)$ if and only if  $\alpha(\mu)\in \Gamma(X)$. Thus,
putting $\zeta_j = \alpha_{1j} \xi_j$, where $j=1, 2, 3$, we
reduce the proof of Theorem 1 to the case when the
 $L_1$, $L_2$, and $L_3$ have the form $L_1 = \xi_1+\xi_2+\xi_3$, $L_2 =
\delta_{21}\xi_1+\delta_{22}\xi_2+\delta_{23}\xi_3$, and $L_3 =
\delta_{31}\xi_1+\delta_{32}\xi_2+\delta_{33}\xi_3$, where
$\delta_{ij} \in {\rm Aut}(X)$. We also note that if $\alpha,
\beta\in{\rm Aut}(X)$, then the linear statistics $L_1$, $L_2$, and
$L_3$ are independent if and only if the linear statistics $L_1$,
$\alpha L_2$, and $\beta L_3$ are independent. Thus, in proving
Theorem 1 we may assume that
 $L_1 =
\xi_1+\xi_2+\xi_3$, $L_2 = \alpha_1\xi_1+\alpha_2\xi_2+\xi_3$, and
$L_3 = \beta_1\xi_1+\beta_2\xi_2+\xi_3$, where  $\alpha_i =
\left(\begin{matrix}a_i&c_i\\ 0&p_i\end{matrix}\right)$, $\beta_i =
\left(\begin{matrix}b_i&d_i\\ 0&q_i\end{matrix}\right),$ and $i=1, 2,$ are
topological automorphisms of the group
$X=\mathbf{R}\times\mathbf{T}$.

Put $\nu_j=\mu_j*\bar{\mu}_j$. Taking into consideration Lemma 7,
we can prove the theorem assuming that all  $\mu_j$ are
nondegenerate distributions. Moreover, it follows from the proof
of Lemma 7, we have in this case  $\widehat\nu_j(s,0)=e^{-\sigma_j
s^2}$,   \ $s\in \mathbf{R},$ where $\sigma_j> 0$, and $j=1, 2,
3$. We note that then $\sigma_j$ satisfy equations
(\ref{t1.2})--(\ref{t1.4}) and statements 1)--5) of Lemma 2 hold.

By Lemma 1 the characteristic functions $\widehat\mu_j(y)$ satisfy
equation (\ref{u1}). Obviously, the characteristic functions
$\widehat\nu_j(y)$ also satisfy equation (\ref{u1}). Put
$\psi_j(y)=-\log \widehat\nu_j(y)$. Then (\ref{u1}) implies that
the functions $\psi_j(y)$ satisfy equation  (\ref{u4}).

By Lemma 6 the function $\psi_1(s, n)$ satisfies equation
(\ref{u13}). Apply now Lemma 7, and note that the inclusion
\begin{equation}\label{x1}
\mathbf{R}\subset N\cap L
\end{equation}
 is always valid. Then it follows from
(\ref{u13}) and (\ref{x1}) that the function $\psi_1(s, n)$
satisfies equation
 (\ref{u16}). Reasoning as above we get that the functions
 $\psi_2(s, n)$ and $\psi_3(s, n)$ also satisfy equation
(\ref{u16}). Applying Lemma 8 to the functions  $\psi_j(s, n)$ we
obtain the representations
\begin{equation}\label{t1.7}
\psi_j(s, n)=\sigma_j s^2+\kappa_j(n)s+\lambda_j(n), \quad
s\in\mathbf{R}, \ n\in\mathbf{Z},
\end{equation}
where $\kappa_j(-n)=-\kappa_j(n)$, $\lambda_j(-n)=\lambda_j(n)$,
 \  $n\in \mathbf{Z}$, and $\kappa_j(0)=0$, $\lambda_j(0)=0, \
j=1, 2, 3$.

By Lemma 6, $\nu=\nu_1*\nu_2*\nu_3\in \Gamma(X)$. Obviously, then
we have
\begin{equation}\label{gf1}
 \widehat\nu(s, n)=\exp\{-(\sigma s^2+\kappa
sn+\lambda n^2)\}, \quad s\in\mathbf{R}, \ n\in\mathbf{Z},
\end{equation}
where $\sigma> 0$, $\kappa\in\mathbf{R}$, $\lambda\geq 0$ and
\begin{equation}\label{t1.7.1}
\psi_1(s, n)+\psi_2(s, n)+\psi_3(s, n)=\sigma s^2+\kappa
sn+\lambda n^2, \quad s\in\mathbf{R}, \ n\in\mathbf{Z}.
\end{equation}

We check that $\kappa_j(n)$ in (\ref{t1.7}) are linear functions.
Put in (\ref{u4}) $u=(0, n)$, $v=(s_2, 0)$, $w=(s_3, 0)$. Note
that $\widetilde\alpha_1 v= (a_1 s_2, 0)$, $\widetilde\alpha_2 v=
(a_2 s_2, 0)$, $\widetilde\beta_1 w= (b_1 s_3, 0)$, and
$\widetilde\beta_2 w= (b_2 s_3, 0)$. Taking into account
representations (\ref{t1.7}), we find from (\ref{u4}) that

\begin{equation}\label{t1.8}
\sigma_1(a_1 s_2+b_1 s_3)^2 + \kappa_1(n)(a_1 s_2+b_1 s_3)
+\lambda_1(n) +
    \sigma_2 (a_2 s_2+b_2 s_3)^2$$ $$ +\kappa_2(n)(a_2 s_2+b_2 s_3)
 +\lambda_2(n)+
    \sigma_3(s_2+s_3)^2 + \kappa_3(n)(s_2+s_3)
+\lambda_3(n)$$
$$ =\lambda_1(n) +\lambda_2(n)+ \lambda_3(n)+\sigma_1(a_1 s_2)^2+\sigma_2
(a_2 s_2)^2+ \sigma_3 s_2^2$$$$+ \sigma_1(b_1 s_3)^2+\sigma_2 (b_2
s_3)^2+ \sigma_3 s_3^2, \quad s_2, s_3\in\mathbf{R}, \
n\in\mathbf{Z}.
\end{equation}
Setting in (\ref{t1.8}) $s_2=1, s_3=0$, we get
\begin{equation}\label{t1.9}
a_1 \kappa_1(n) +
    a_2 \kappa_2(n) + \kappa_3(n) =0, \
n\in\mathbf{Z}.
\end{equation}
Setting in (\ref{t1.8}) $s_2=0, s_3=1$, we find
\begin{equation}\label{t1.10}
b_1 \kappa_1(n) +
    b_2 \kappa_2(n) + \kappa_3(n) =0, \quad
n\in\mathbf{Z}.
\end{equation}
On the other hand it follows from   (\ref{t1.7}) and
(\ref{t1.7.1}) that
\begin{equation}\label{t1.11}
\kappa_1(n) + \kappa_2(n) + \kappa_3(n) =\kappa n, \quad
n\in\mathbf{Z}.
\end{equation}
Subtracting from  (\ref{t1.9}) and (\ref{t1.10})
equality (\ref{t1.11}), we obtain
\begin{equation}\label{t1.12}
(a_1-1) \kappa_1(n) +
    (a_2-1) \kappa_2(n)=-\kappa n, \quad
n\in\mathbf{Z}.
\end{equation}
\begin{equation}\label{t1.13}
(b_1-1) \kappa_1(n) +
    (b_2-1) \kappa_2(n)=-\kappa n, \quad
n\in\mathbf{Z}.
\end{equation}
Taking into account that by statement  5) of Lemma 2 the
inequality $\left|
\begin{array}{cc}
  a_1-1 & a_2-1 \\
  b_1-1 & b_2-1 \\
\end{array}
\right|\neq 0$ holds, it follows from (\ref{t1.11})--(\ref{t1.13})
that $\kappa_j(n)$ are linear functions, i.e.
$\kappa_j(n)=\kappa_j n, \ n\in\mathbf{Z}$, where $\kappa_j\in
\mathbf{R}$. Thus, representations (\ref{t1.7}) take the form
\begin{equation}\label{t1.14}
\psi_j(s, n)=\sigma_j s^2+\kappa_j sn+\lambda_j(n), \quad
s\in\mathbf{R}, \ n\in\mathbf{Z}, \ j=1, 2, 3.
\end{equation}
 Note that it follows from  (\ref{t1.7.1}) and (\ref{t1.14})  that
\begin{equation}\label{t1.15}
\lambda_1(n)+\lambda_2(n)+\lambda_3(n)=\lambda n^2, \quad
n\in\mathbf{Z}.
\end{equation}

Let $u=(s_1, n_1)$, $v=(s_2, n_2)$, $w=(s_3, n_3)$. Then
$\widetilde\alpha_1 v= (a_1 s_2+ c_1 n_2, p_1 n_2)$,
$\widetilde\alpha_2 v= (a_2 s_2+ c_2 n_2, p_2 n_2)$,
$\widetilde\beta_1 w= (b_1 s_3+ d_1 n_3, q_1 n_3)$,
$\widetilde\beta_2 w= (b_2 s_3+ d_2 n_3, q_2 n_3)$. Substituting
representations (\ref{t1.14}) into (\ref{u4}), we get that
equalities (\ref{t1.2})--(\ref{t1.4}) hold. Moreover, the
following equalities also hold:
\begin{equation}\label{t1.17}
    \kappa_1a_1+\kappa_2a_2+\kappa_3=0
\end{equation}
\begin{equation}\label{t1.18}
    \kappa_1b_1+\kappa_2b_2+\kappa_3=0
\end{equation}
\begin{equation}\label{t1.19}
    2\sigma_1c_1 +2\sigma_2c_2 + \kappa_1p_1 +\kappa_2 p_2 +\kappa_3=0
\end{equation}
\begin{equation}\label{t1.20}
    2\sigma_1d_1 + 2\sigma_2d_2 + \kappa_1q_1+\kappa_2q_2+\kappa_3=0
\end{equation}
\begin{equation}\label{t1.21}
    2\sigma_1a_1d_1 + 2\sigma_2a_2d_2 +
\kappa_1a_1q_1+\kappa_2a_2q_2+\kappa_3=0
\end{equation}
\begin{equation}\label{t1.22}
    2\sigma_1b_1c_1 + 2\sigma_2b_2c_2 +
 \kappa_1b_1p_1+\kappa_2b_2p_2+\kappa_3=0
\end{equation}
\begin{equation}\label{t1.23}
    n_1n_2(\kappa_1c_1+\kappa_2c_2)+n_1n_3(\kappa_1d_1+\kappa_2d_2)
    $$$$+n_2n_3(2\sigma_1c_1d_1+2\sigma_2c_2d_2+\kappa_1c_1q_1+\kappa_1d_1p_1
+\kappa_2c_2q_2+\kappa_2d_2p_2)$$ $$+\lambda_1(n_1+p_1 n_2+q_1
n_3)+\lambda_2(n_1+p_2 n_2+q_2 n_3)$$$$+\lambda_3(n_1+n_2+n_3) =
\lambda( n_1^2+n_2^2+n_3^2), \quad n_j\in\mathbf{Z}.
\end{equation}

We find from (\ref{t1.7.1})  and (\ref{t1.14}) that
\begin{equation}\label{f7}
\sigma_1+\sigma_2+\sigma_3=\sigma,
\end{equation}
and
\begin{equation}\label{f8}
\kappa_1+\kappa_2+\kappa_3=\kappa.
\end{equation}

Prove that the support of $\nu$ is a subgroup of  $X$
topologically isomorphic to  $ \mathbf{R}$. It follows from
(\ref{gf1}) that this statement will be proved if we check that
\begin{equation}\label{f11}
4\sigma\lambda=\kappa^2,
\end{equation}
i.e. taking into account (\ref{f7})  and (\ref{f8}), we should
prove that
\begin{equation}\label{t1.24}
    4(\sigma_1+\sigma_2+\sigma_3)\lambda=(\kappa_1+\kappa_2+\kappa_3)^2.
\end{equation}
The proof of (\ref{t1.24}) is the series of elementary
and boring computations. Note that by statement 4) of Lemma 2 we
have ${a_2b_1-a_1b_2\neq 0}$. Then it follows from   (\ref{t1.2})
and (\ref{t1.3}) that
\begin{equation}\label{t1.25}
\sigma_1={b_2-a_2\over a_2b_1-a_1b_2} \sigma_3,\quad
\sigma_2={a_1-b_1\over a_2b_1-a_1b_2} \sigma_3.
\end{equation}
Similarly we find from (\ref{t1.17}) and (\ref{t1.18})
that
\begin{equation}\label{t1.26}
\kappa_1={b_2-a_2\over a_2b_1-a_1b_2} \kappa_3,\quad
\kappa_2={a_1-b_1\over a_2b_1-a_1b_2} \kappa_3.
\end{equation}

Each of the numbers  $p_j$ and $q_j$ can take the values $\pm 1$.
We have here some cases and consider  each of them separately.

\textbf{I.} Let $p_1=p_2=1$. Then (\ref{t1.19}) implies that
\begin{equation}\label{t1.27}
    \kappa_1+\kappa_2+\kappa_3=-2\sigma_1c_1 -2\sigma_2c_2
\end{equation}
Putting in (\ref{t1.23}) $n_1=-n_2=1, \ n_3=0$, we get
\begin{equation}\label{t1.28}
    -\kappa_1c_1-\kappa_2c_2=2\lambda.
\end{equation}
Substituting  (\ref{t1.27}) and (\ref{t1.28}) into
(\ref{t1.24}) and taking into account (\ref{t1.25}) and
(\ref{t1.26}), we make sure that we obtain an equality. Thus, in
case \textbf{I} equality (\ref{f11}) is proved.

\textbf{II.} Let $q_1=q_2=1$. Then (\ref{t1.20}) implies that
\begin{equation}\label{t1.31}
    \kappa_1+\kappa_2+\kappa_3=-2\sigma_1d_1 - 2\sigma_2d_2
\end{equation}
Putting in (\ref{t1.23}) $n_1=-n_3=1, n_2=0$, we obtain
\begin{equation}\label{t1.32}
    -\kappa_1d_1-\kappa_2d_2=2\lambda.
\end{equation}
Substituting (\ref{t1.31}) and (\ref{t1.32}) into
(\ref{t1.24}) and taking into account (\ref{t1.25}) and
(\ref{t1.26}), we make sure that we obtain an equality.
 Thus, in case \textbf{II} equality (\ref{f11}) is proved.

\textbf{III.} Let either $p_1=p_2=q_1=q_2=-1$  or $p_1=q_1=-1$ and
$p_2=q_2=1$  or $p_1=q_1=1$ and $p_2=q_2=-1$. Putting in
(\ref{t1.23}) $n_1=0, \ n_2=-n_3=1$, we get
\begin{equation}\label{t1.35}
    -(2\sigma_1c_1d_1+2\sigma_2c_2d_2+\kappa_1c_1p_1+\kappa_1d_1p_1
+\kappa_2c_2p_2+\kappa_2d_2p_2)=2\lambda.
\end{equation}
Note that by statement 3) of Lemma 2 we have $a_1\ne
a_2$ and $b_1\ne b_2$. It follows from   (\ref{t1.19}) and
(\ref{t1.22}) that
\begin{equation}\label{t1.35.1}
    c_1= {(1-b_2)\kappa_3- \kappa_1p_1 (b_2-b_1)\over
    2\sigma_1(b_2-b_1)}, \quad c_2= {(1-b_1)\kappa_3- \kappa_2p_2 (b_1-b_2)\over
    2\sigma_2(b_1-b_2)}.
\end{equation}
We find from (\ref{t1.20}) and (\ref{t1.21}) that
\begin{equation}\label{t1.35.2}
    d_1= {(1-a_2)\kappa_3- \kappa_1q_1 (a_2-a_1)\over  2\sigma_1(a_2-a_1)},
    \quad d_2= {(1-a_1)\kappa_3- \kappa_2q_2 (a_1-a_2)\over 2\sigma_2(a_1-a_2)}.
\end{equation}
Finding $\lambda$ from  (\ref{t1.35})  and substituting
the obtained expression into (\ref{t1.24}), we get that we should
check the equality
\begin{equation}\label{t1.35.3}
    -2(2\sigma_1c_1d_1+2\sigma_2c_2d_2+\kappa_1c_1p_1+\kappa_1d_1p_1
+\kappa_2c_2p_2+\kappa_2d_2p_2)$$$$\times(\sigma_1+\sigma_2+\sigma_3)=(\kappa_1+\kappa_2+\kappa_3)^2.
\end{equation}
Substitute (\ref{t1.25}), (\ref{t1.26}),
(\ref{t1.35.1}), and (\ref{t1.35.2}) into (\ref{t1.35.3}). After
elementary computations we see that the verification of
(\ref{t1.35.3}) is reduced to the verification of the equality
\begin{equation}\label{t1.35.3+2}
    (a_2b_1-a_1b_2)(
    (1-b_2)(1-a_2)(a_1-b_1)+(1-b_1)(1-a_1)(b_2-a_2)
    )$$$$=-(b_2-b_1)(a_2-a_1)(a_1-b_1)(b_2-a_2).
\end{equation}
The correctness of (\ref{t1.35.3+2}) follows from statement 1) of
Lemma 2. Thus, in case \textbf{III} equality (\ref{f11}) is
proved.

To consider the remaining cases  put in (\ref{t1.23})
$n_1=n_2=n_3=1$. Taking into account that
$\lambda_j(-n)=\lambda_j(n)$,
 \  $n\in \mathbf{Z}$, \ $
j=1, 2$, we get
\begin{equation}\label{t1.23.1}
    \kappa_1c_1+\kappa_2c_2+\kappa_1d_1+\kappa_2d_2+2\sigma_1c_1d_1+2\sigma_2c_2d_2+\kappa_1c_1q_1$$
$$+\kappa_1d_1p_1
+\kappa_2c_2q_2+\kappa_2d_2p_2+\lambda_1(1)+\lambda_2(1)+\lambda_3(3)
= 3\lambda.
\end{equation}

\textbf{IV.} Let $p_1=p_2=-1$,  and $q_1=-q_2$. Putting in
(\ref{t1.23}) $n_1=2, n_2=1, n_3=0$, we obtain
$$
2(\kappa_1c_1+\kappa_2c_2)+\lambda_1(1)+\lambda_2(1)+\lambda_3(3) = 5\lambda.
$$
It follows from this and from (\ref{t1.23.1}) that
$$
\kappa_1c_1+\kappa_2c_2-(2\sigma_1c_1d_1+2\sigma_2c_2d_2+\kappa_1c_1q_1
-\kappa_2c_2q_1)= 2\lambda.
$$

\textbf{V.} Let $q_1=q_2=-1$,  and $p_1=-p_2$. Putting in
(\ref{t1.23}) $n_1=2, n_2=0, n_3=1$, we get

$$
2(\kappa_1d_1+\kappa_2d_2)+\lambda_1(1)+\lambda_2(1)+\lambda_3(3) = 5\lambda.
$$
It follows from this and from (\ref{t1.23.1}) that
$$
\kappa_1d_1+\kappa_2d_2-(2\sigma_1c_1d_1+2\sigma_2c_2d_2+\kappa_1d_1p_1
-\kappa_2d_2p_1)=2\lambda.
$$

\textbf{VI.} Let either $p_1=q_2=-1$, $p_2=q_1=1$, or $p_1=q_2=1$,
$p_2=q_1=-1$.
 Putting in (\ref{t1.23}) $n_1=0, n_2=2, n_3=1$, we find
$$
2(2\sigma_1c_1d_1+2\sigma_2c_2d_2+\kappa_1c_1q_1+\kappa_1d_1p_1
+\kappa_2c_2q_2+\kappa_2d_2p_2)$$$$+\lambda_1(1)+\lambda_2(1)+\lambda_3(3)
= 5\lambda.
$$
It follows from this and from (\ref{t1.23.1}) that
$$
2\sigma_1c_1d_1+2\sigma_2c_2d_2+\kappa_1c_1p_2+\kappa_1d_1p_1
+\kappa_2c_2p_1+\kappa_2d_2p_2$$$$-(\kappa_1d_1+\kappa_2d_2)-(\kappa_1c_1+\kappa_2c_2)
=2\lambda.
$$
Then in each of these cases we reason as in case
 \textbf{III.} Thus, we have proved that equality  (\ref{f11}) holds.
 So, we have obtained the representation
$$
\widehat\nu(s, n)=\exp\{-(\sigma s^2+\kappa sn+\lambda
n^2)\}=\exp\left\{ -\sigma\left(s+{\kappa\over 2\sigma}
n\right)^2\right\}, \quad s\in\mathbf{R}, \ n\in\mathbf{Z}.
$$

Set $H=\{y\in Y:\ \widehat\nu(y)=1\}$. It is obvious that
$H=\{n(-{\kappa\over 2\sigma},1):\ n\in \mathbf{Z} \}$. Then
$G=A(X,H)=\{(t, e^{it{\kappa\over 2\sigma}}): \ t\in\mathbf{R}\}$.
It is clear that $G\cong \mathbf{R}$. We have $\sigma(\nu)\subset
G$. Since
$\nu=\nu_1*\nu_2*\nu_3=\mu_1*\bar{\mu}_1*\mu_2*\bar{\mu}_2*\mu_3*\bar{\mu}_3$,
it follows from Lemma 3 that   the distributions  $\mu_j$ may be
replaced by their shifts
 ${\mu}'_j$ such that
$\nu=\mu'_1*{\bar\mu}'_1*\mu'_2*{\bar\mu}'_2*\mu'_3*{\bar\mu}'_3$,
and $\sigma(\mu'_j)\subset G$. Since $\nu$ is a Gaussian
distribution and  $G\cong \mathbf{R}$, by the classical Cram\'er
theorem on the decomposition of a Gaussian distribution, we get
that  all $\mu'_j$ Gaussian distributions, and hence  $\mu_j$ are
also Gaussian. $\Box$

\bigskip

{\textbf{5. Comments to Theorem 1.}

 \textbf{Remark 1.}  For $n=2$ on the group
$X=\mathbf{R}\times\mathbf{T}$ under the conditions of Problem 1
 the
distributions $\mu_j$ are
either Gaussian or convolutions of Gaussian
distributions and signed measures supported in  $\mathbf{Z}(2)$
\cite{FeMy1}. According to Theorem 1 for $n=3$ the
distributions $\mu_j$ are only Gaussian. We see that in contrast
to the case of the real line $\mathbf{R}$ on the group
$X=\mathbf{R}\times\mathbf{T}$ the class of distributions which are
characterized by the independence of  $n$ linear statistics of $n$
independent random variables depends on $n$.

 \textbf{Remark 2.} As   established in the proof of Theorem 1,
 if   $L_1$, $L_2$, and $L_3$  are of the form $L_1 = \xi_1+\xi_2+\xi_3$,
$L_2 = \alpha_1\xi_1+\alpha_2\xi_2+\xi_3$, $L_3 =
\beta_1\xi_1+\beta_2\xi_2+\xi_3$, and $\mu_j$ are nondegenerate
distributions, then supports of  $\mu_j$ are cosets of the group $X$ with respect
the one-parameter subgroup  $G={\{(t, e^{it{\kappa\over 2\sigma}}): \
t\in\mathbf{R}\}}$.

We will check that the subgroup  $G$ is invariant with respect to
all topological automorphisms   $\alpha_i, \beta_i$. Taking into
account (\ref{f7}), (\ref{f8}), (\ref{t1.25}) and (\ref{t1.26}),
we get that  $G=\{(t,e^{it{\kappa_3\over 2\sigma_3}}): \
t\in\mathbf{R}\}$. Verify that the subgroup  $G$ is invariant with
respect to $\alpha_1$. Since
$$
\alpha_1(t,e^{it{\kappa_3\over 2\sigma_3}})=
\left(a_1 t, e^{ic_1 t} e^{itp_1{\kappa_3\over 2\sigma_3}}\right)=
\left(a_1 t, e^{ia_1 t {\kappa_3\over 2\sigma_3}
\left({2c_1\sigma_3\over a_1\kappa_3} +{p_1\over a_1}\right) }
\right),
$$
it suffices to show that
\begin{equation}\label{f9}
    {2c_1\sigma_3\over a_1\kappa_3} +{p_1\over a_1}=1.
\end{equation}
To check equality (\ref{f9}) substitute first in (\ref{f9}) the
representation for $c_1$ from (\ref{t1.35.1}), and then substitute
the representation for $\sigma_1$ and $\kappa_1$ from
(\ref{t1.25}) and (\ref{t1.26}) into the obtained expression.
After simple transformations we see that (\ref{f9}) is reduced to
equality  1) of Lemma 2.

Reasoning as above we get that the subgroup  $G$ is invariant with
respect to    $\alpha_2, \beta_1, \beta_2$. Obviously, the
restriction of any automorphism   $\alpha_i, \beta_i$ to the
subgroup $G$ is a topological automorphism of  $G$.

It follows from what has been said that the statement   in Theorem
1 that $\mu_j$
 are Gaussian distributions can be made without the
Cram\'er theorem but directly from the Skitovich--Darmois theorem
for the real line.

We note that the  property:  $G$ is invariant with respect to all
topological automorphisms   $\alpha_i, \beta_i$, imposes the
strong restrictions on $\alpha_i, \beta_i$. Indeed, let
  $G$ be an arbitrary subgroup of the group
$X=\mathbf{R}\times\mathbf{T}$ topologically isomorphic to
$\mathbf{R}$. Then $G$ is of the form $G=\{(t,e^{i\omega t}): \
t\in\mathbf{R}\}$, where $\omega$ is a fixed real number. It is
obvious that the subgroup $G$ is invariant with respect to
an automorphism  $\alpha=\left(%
\begin{array}{cc}
  a & c \\
  0 & p \\
\end{array}%
\right) \in \mathrm{Aut} (X)$ if and only if $c=(a-p)\omega$.

Let  $\omega=0$, i.e. $ G=\{(t, 1): \ t\in\mathbf{R}\}$. Then
$c_1=c_2=d_1=d_2=0$, i.e.
  $\alpha_i, \beta_i$ are diagonal matrices, moreover statements 1)--5) of Lemma 2 hold
  for  $a_i$ and $b_i$.

 Let  $\omega\ne 0$. Assume also that no
  $\alpha_i, \beta_i$ is equal to $\pm I$. It follows from this that
  in addition to statements
1)--5) of Lemma  2 the equalities
\begin{equation}\label{f10}
    {c_1\over a_1- p_1}={c_2\over a_2- p_2}={d_1\over b_1- q_1}={d_2\over b_2-
    q_2}
\end{equation}
also hold. In particular, if at least one of these equalities does
not hold, then all   $\mu_j$ are degenerate distributions. If some
of the automorphisms $\alpha_i, \beta_i$ are equal to $\pm I$,
then in (\ref{f10}) the expressions corresponding to these
automorphisms are omitted.

\textbf{Remark 3.} Theorem 1 can not strengthened to the statement
that  $\mu_j$ are degenerate distributions (compare with Lemma 5).
Indeed, let $G$ be an arbitrary one-parameter subgroup of the
group  $X=\mathbf{R}\times\mathbf{T}$ of the form
$G=\{(t,e^{i\omega t}): \ t\in\mathbf{R}\}$, where $\omega$ is a
fixed real number. Let $a_i, b_i$, where $i = 1, 2,$ be nonzero
real numbers such that system of equations
(\ref{t1.2})--(\ref{t1.4}) has a solution $\sigma_1>0, \sigma_2>0,
\sigma_3>0$. Consider the automorphisms  $\alpha_i, \beta_i \in
\mathrm{Aut}(X)$ of the form
$$\alpha_1=\left(\begin{matrix}a_1& (a_1-p_1)\omega\\ 0&p_1\end{matrix}\right),
\quad \alpha_2=\left(\begin{matrix}a_2&(a_2-p_2)\omega\\
0&p_2\end{matrix}\right), \
$$$$\beta_1 = \left(\begin{matrix}b_1&(b_1-q_1)\omega\\ 0&q_1\end{matrix}\right),
\quad \beta_2 = \left(\begin{matrix}b_2&(b_2-q_2)\omega\\
0&q_2\end{matrix}\right).$$ Let $\xi_j$, where $j=1, 2, 3,$ be independent
Gaussian random variables with values in the group $X$ and
distributions  $\mu_j$ having the characteristic functions
$$ \widehat\mu_1(s,
n)=\exp\{-\sigma_1(s+\omega n)^2\}, \quad \widehat\mu_2(s,
n)=\exp\{- \sigma_2(s+\omega n)^2\},
$$$$\widehat\mu_3(s, n)=\exp\{-\sigma_3(s+\omega n)^2\}, \quad s\in\mathbf{R}, \ n\in\mathbf{Z}.$$
Obviously, the support  of any distribution  $\mu_j$ coincides
with the subgroup  $G$. Taking into account
(\ref{t1.2})--(\ref{t1.4}), it is easy to verify directly that the
characteristic functions   $ \widehat\mu_j(s, n)$ satisfy equation
  (\ref{u1}), and hence, by Lemma 1 the linear statistics
$L_1 = \xi_1+\xi_2+\xi_3$, $L_2 =
\alpha_1\xi_1+\alpha_2\xi_2+\xi_3$, and $L_3 =
\beta_1\xi_1+\beta_2\xi_2+\xi_3$ are independent.

\textbf{Remark 4.} On the group $X=\mathbf{T}$ for $n=2$ under the
conditions of Problem 1  the distributions $\mu_j$
are
either Gaussian or convolutions of Gaussian
distributions and signed measures
\cite{Ba-E-Sta}, and for $n=3$ all
$\mu_j$ are degenerate distributions by Lemma 5. On the group
${X=\mathbf{R}\times\mathbf{T}}$ for $n=2$ under the conditions of
Problem 1  the distributions $\mu_j$ are also
either Gaussian or convolutions of Gaussian
distributions and signed measures \cite{FeMy1}, but
for $n=3$ all $\mu_j$ are Gaussian distributions by Theorem 1. So, we could
surmise that with increasing $n$, the class of distributions in
Problem 1, which are characterized by the independence of  $n$
linear statistics of $n$ independent random variables, decreases. It
turns out that this assumption is not true.

Indeed, let $\xi_1$, $\xi_2$, $\xi_3$, $\xi_4$ be independent
random variables with values in the group $\mathbf{T}$ and
distributions  $\mu_j$. Consider the linear statistics
$L_1=\xi_1+\xi_2+\xi_3+\xi_4$, $L_2=\xi_1+\xi_2-\xi_3-\xi_4$,
$L_3=\xi_1-\xi_2+\xi_3-\xi_4$, and $L_4=\xi_1-\xi_2-\xi_3+\xi_4$.
Take $\mu_1=\mu_2=\gamma*\pi_1$, $\mu_3=\mu_4=\gamma*\pi_2$, where
$\gamma\in \Gamma(\mathbf{T})$, and $\pi_j$ are signed measures
supported in  ${\mathbf Z}(2)$ such that $\pi_1*\pi_2=E_1$. It is
easy to verify that in this case the linear statistics  $L_1$, $L_2$,
$L_3$, and $L_4$ are independent, but all
$\mu_j\not\in\Gamma(\mathbf{T})$.

\bigskip

{\textbf{6. Solution of Problem  1 for the group
$X=\Sigma_\text{\boldmath $a$}\times\mathbf{T}.$}} Pass    to the
solution of Problem 1  in the case  $n=3$ for the group
${X=\Sigma_\text{\boldmath $a$}\times\mathbf{T}}$. Remind the
definition of $\text{\boldmath $a$}$-adic solenoid
$\Sigma_\text{\boldmath $a$}$.

Let $\text{\boldmath $a$}= (a_0, a_1,\dots)$, where all $a_j \in {\mathbf{Z}}$,
$a_j > 1$. First we define the   group of $\text{\boldmath $a$}$-adic
integers $\Delta_\text{\boldmath $a$}$. As a set $\Delta_\text{\boldmath $a$}$ coincides
with the Cartesian product $\mathop{\mbox{\rm\bf
P}}\limits_{n=0}^\infty\{0,1,\dots ,a_n-1\}$. Consider
$\text{\boldmath $x$}=(x_0, x_1, x_2,\dots)$, $\text{\boldmath $y$}=(y_0, y_1,
y_2,\dots)\in \Delta_\text{\boldmath $a$}$, and define the sum
$\text{\boldmath $z$}=\text{\boldmath $x$}+\text{\boldmath $y$}$ as follows. Let
$x_0+y_0=t_0a_0+z_0$, where $z_0\in\{0, 1,\dots ,a_0-1\}$, $
t_0\in\{0, 1\}$. Assume that the numbers $z_0,z_1,\dots ,z_k;$
$t_0, t_1,\dots ,t_k$ have been already determined. Let us put then
$x_{k+1}+y_{k+1}+t_k=t_{k+1}a_{k+1}+z_{k+1}$, where
$z_{k+1}\in\{0,1,\dots ,a_{k+1}-1\}$, $t_{k+1}\in\{0, 1\}$. This
defines by induction a sequence $\text{\boldmath $z$}=(z_0, z_1,
z_2,\dots)$. The set $\Delta_\text{\boldmath $a$}$ with the addition defined
above is an Abelian group.  The obtained group considering in the product topology
 is called the $\text{\boldmath $a$}$-adic integers. Consider the group
$\mathbf{R}\times\Delta_\text{\boldmath $a$}$. Let $B$ be the subgroup of the
group $\mathbf{R}\times\Delta_\text{\boldmath $a$}$ of the form
$B=\{(n,n\text{\boldmath $u$})\}_{n=-\infty}^{\infty}$, where
$\text{\boldmath $u$}=(1, 0,\dots,0,\dots)$. The factor-group
$\Sigma_\text{\boldmath $a$}=(\mathbf{R}\times\Delta_\text{\boldmath $a$})/B$ is called
the $\text{\boldmath $a$}$-{adic solenoid}.  The group $\Sigma_\text{\boldmath $a$}$ is compact,
 connected and has dimension one \cite[(10.12), (10.13),
(24.28)]{Hewitt-Ross}.
Denote by
$\mathbf{Q}$ the group of rational numbers considering in the discrete topology.
The character group of the group $\Sigma_\text{\boldmath $a$}$ is topologically
isomorphic to a subgroup $H_\text{\boldmath $a$}\subset
\mathbf{Q}$ of the form
$$H_\text{\boldmath $a$}= \left\{{m \over a_0a_1 \dots a_n} : \quad n = 0, 1,\dots; \ m
\in {\mathbb{Z}} \right\}
$$
\cite[(25.3)]{Hewitt-Ross}.

Let $X=\Sigma_\text{\boldmath $a$}\times\mathbf{T}$. Denote by $x=(g,z)$,
where $g\in\Sigma_\text{\boldmath $a$}$, and $z\in\mathbf{T}$, elements of
the group $X$. If $X=\Sigma_\text{\boldmath $a$}\times\mathbf{T}$,
 then $Y\cong H_\text{\boldmath $a$}\times\mathbf{Z}$.  In order
to avoid introducing new notation we assume that
$Y=H_\text{\boldmath $a$}\times\mathbf{Z}$,   and denote by $y=(r,n)$, where
$r\in H_\text{\boldmath $a$}$, and $n\in\mathbf{Z}$, elements of $Y$. It is
easy to verify that every automorphism $\varepsilon\in{\rm
Aut}(H_\text{\boldmath $a$}\times\mathbf{Z})$ is defined by a matrix
 $\left(\begin{matrix}a &c\\ 0&p\end{matrix}\right)$, where $a\in {\rm
Aut}(H_\text{\boldmath $a$}), \ c\in H_\text{\boldmath $a$}, p=\pm 1,$  and
$\varepsilon$ operates on $Y=H_\text{\boldmath $a$}\times\mathbf{Z}$ as
follows
$$\varepsilon(r,n)=(a r+c n, pn),
\quad r\in H_\text{\boldmath $a$}, \ n\in\mathbf{Z}.
$$
Then the
adjoint automorphism $\delta=\widetilde{\varepsilon}\in{\rm
Aut}(X)$ is of the form
$$ \delta(g,z)=(\widetilde{a} g, (g,c)z^{p}),
\quad g\in\Sigma_\text{\boldmath $a$}, \ z\in\mathbf{T}.
$$

 Denote by
  $\iota$  the natural embedding $\iota :Y\mapsto
\mathbf{R}\times\mathbf{Z}, \ \iota (r, n)=(r, n).$ Let
$\tau=\widetilde\iota$    be the adjoint homomorphism
$\tau:\mathbf{R}\times\mathbf{T}\mapsto X$. Put $g_t=\tau (t, 1),
\ g_t\in \Sigma_\text{\boldmath $a$}.$ Then $\tau(t, z)=(g_t, z)$. Since,
 $\overline{\iota (Y)}=\mathbf{R}\times\mathbf{Z}$, we conclude
 that
$\tau$
 is a monomorphism  \cite[(24.41)]{Hewitt-Ross}.

\textbf{Lemma 9} (\cite{FeMy1}). \textit{Let $r$ be a rational
number, $\tau:\mathbf{R}\times\mathbf{T}\mapsto X$ be the
homomorphism defined above. Let $K=\{(t, \ e^{itr}): \
t\in\mathbf{R}\}$ be the subgroup of $\mathbf{R}\times\mathbf{T}$.
Then $\overline{\tau(K)}\cong\Sigma_\text{\boldmath $a$}$.}

\textbf{Lemma 10}.  \textit{Let
$X=\Sigma_\text{\boldmath $a$}\times\mathbf{T}$, and $\alpha_{ij}\in
\mathrm{Aut}(X)$, where $i, j=1,2,3$. Let $\xi_j$, where $j=1, 2,
3$, be independent random variables with values in $X$ and
distributions $\mu_j$ such that their characteristic functions do
not vanish. Assume that the linear statistics
$L_1=\alpha_{11}\xi_1+\alpha_{12}\xi_2+\alpha_{13}\xi_3$,
$L_2=\alpha_{21}\xi_1+\alpha_{22}\xi_2+\alpha_{23}\xi_3$, and
$L_3=\alpha_{31}\xi_1+\alpha_{32}\xi_2+\alpha_{33}\xi_3$ are
independent.  Then the distributions $\mu_j$ may be replaced by
their shifts  ${\mu}'_j$ such that ${\mu}'_j$ are concentrated on
the subgroup $\tau(\mathbf{R}\times\mathbf{T})$.}

The proof of this lemma can be done in a way as in
\cite{FeMy1} the corresponding statement was proved for two
independent random variables $\xi_1$ and $\xi_2$ with values in
the group $X$ and two linear statistics $L_1=\xi_1+\xi_2$ and
$L_2=\xi_1+\delta\xi_2$, where ${\delta}\in {\rm Aut}(X)$. $\Box$

\textbf{Theorem 2.} \textit{Let
$X=\Sigma_\text{\boldmath $a$}\times\mathbf{T}$, and $\alpha_{ij}\in
\mathrm{Aut}(X)$, where $i, j=1,2,3$. Let $\xi_j$, where $j=1, 2,
3$, be independent random variables with values in $X$ and
distributions $\mu_j$ such that their characteristic functions do
not vanish. Assume that the linear statistics
$L_1=\alpha_{11}\xi_1+\alpha_{12}\xi_2+\alpha_{13}\xi_3$,
$L_2=\alpha_{21}\xi_1+\alpha_{22}\xi_2+\alpha_{23}\xi_3$, and
$L_3=\alpha_{31}\xi_1+\alpha_{32}\xi_2+\alpha_{33}\xi_3$ are
independent. Then  all $\mu_j$ are either degenerate distributions
or
 Gaussian distributions such that their supports are cosets of
 a subgroup of $X$ topologically
isomorphic to    $\Sigma_\text{\boldmath $a$}$.}

\textbf{Proof.} Let the homomorphism $\tau$ be defined as above.
Put $B=\tau(\mathbf{R}\times \mathbf{T})$. It is easily seen that
the subgroup $B$ is invariant with respect to any automorphism
$\alpha\in{\rm Aut}(X)$. For each automorphism $\alpha$ we define
the mapping   $\bar\alpha : \mathbf{R}\times \mathbf{T}\mapsto
\mathbf{R}\times \mathbf{T}$ by the formula $\bar\alpha (t,
z)=\tau^{-1}\alpha\tau (t, z)$. It is easily verified that
$\bar\alpha \in{\rm Aut}(\mathbf{R}\times \mathbf{T})$, and the
same matrix corresponds to $\bar\alpha$ and to $\alpha$. Taking
into account Lemma 10, we can assume from the beginning that the
distributions $\mu_j$ are concentrated on the subgroup  $B$. Note
that because
 $\tau$ is a continuous and on-to-one mapping, by the Suslin theorem
  images of Borel sets under the mapping
$\tau$ are also Borel sets.

Let $\xi$ be a random variable  with values in the group  $X$ and
with the distribution  concentrated on the subgroup  $B$. Put
$\widetilde\xi=\tau^{-1}\xi$. Then $\widetilde\xi$ is a random
variable  with values in the group $\mathbf{R}\times\mathbf{T}$.
Consider random variables $\widetilde\xi_j$. Since
$\tau^{-1}\alpha_{ij}\xi_j=\bar\alpha_{ij}\widetilde\xi_j$, the
linear statistics $\widetilde
L_1=\bar\alpha_{11}\widetilde\xi_1+\bar\alpha_{12}\widetilde\xi_2+\bar\alpha_{13}\widetilde\xi_3$,
$\widetilde L_2=\bar\alpha_{21}\widetilde
\xi_1+\bar\alpha_{22}\widetilde\xi_2+\bar\alpha_{23}\widetilde\xi_3$,
and $\widetilde L_3=\bar\alpha_{31}\widetilde
\xi_1+\bar\alpha_{32}\widetilde\xi_2+\bar\alpha_{33}\widetilde\xi_3$
are independent. Let $\widetilde\mu_j$ be the distributions of the
random variables $\widetilde\xi_j$. Since
$\mu_j=\tau(\widetilde\mu_j)$,  Theorem 2 follows from Theorem 1
and Lemma 9. $\Box$

\vskip 1 cm

\noindent G.M. Feldman\\
Mathematical Division,\\ B. Verkin Institute for Low Temperature
Physics
 and Engineering \\of the National Academy of Sciences of
 Ukraine,\\
47, Lenin ave, Kharkov, 61103, Ukraine

\bigskip

\noindent e-mail:
feldman@ilt.kharkov.ua 

\vskip 1 cm

\noindent M.V. Myronyuk\\
Mathematical Division,\\ B. Verkin Institute for Low Temperature
Physics
 and Engineering \\of the National Academy of Sciences of
 Ukraine,\\
47, Lenin ave, Kharkov, 61103, Ukraine

\bigskip

\noindent e-mail:
myronyuk@ilt.kharkov.ua 
\end{document}